\documentclass{article}

\usepackage{PRIMEarxiv}

\usepackage[utf8]{inputenc} 
\usepackage[T1]{fontenc}    
\usepackage{hyperref}       
\usepackage{url}            
\usepackage{booktabs}       
\usepackage{amsfonts}       
\usepackage{nicefrac}       
\usepackage{microtype}      
\usepackage{lipsum}
\usepackage{fancyhdr}       
\usepackage{graphicx}       
\graphicspath{{media/}}     
\usepackage{amsmath} 
\title{Hermite Neural Network Simulation for Solving the 2D Schrodinger Equation

}

\author{
  Kourosh Parand *, Aida Pakniyat \\
  Department of Computer and Data Sciences, Faculty of Mathematical Sciences,\\
  Shahid Beheshti University \\
  Tehran, Iran\\
  {Corresponding author*: Kourosh Parand}
  \texttt{\{Kourosh Parand\}k\_parand@sbu.ac.ir} \\
}

\begin{document}
\maketitle

\begin{abstract}

The Schrodinger equation is a mathematical equation describing the wave function's behavior in a quantum-mechanical system.  It is a partial differential equation that provides valuable insights into the fundamental principles of quantum mechanics. 
In this paper, the aim was to solve the Schrodinger equation with sufficient accuracy by using a mixture of neural networks with the collocation method base Hermite functions. Initially, the Hermite functions roots were employed as collocation points, enhancing the efficiency of the solution. The Schrodinger equation is defined in an infinite domain, the use of Hermite functions as activation functions resulted in excellent precision. Finally, the proposed method was simulated using MATLAB's Simulink tool. 
The results were then compared with those obtained using Physics-informed neural networks and the presented method.
\end{abstract}


\section{Introduction}
Differential equations are used in various scientific and engineering fields, and different methods, such as numerical methods, are applied to solve them.
Quantum mechanics is everywhere, and it has played a fundamental role in developing our understanding of the universe. Many equations are involved in this field, the primary goal being to accurately describe objects' physical properties at the atomic and subatomic levels using the Schrodinger equation. Before this, it was expected that particles of roughly the same size would behave similarly to objects in classical mechanics. However, to everyone's surprise, not only was this not true, but the reality was far stranger than anything that could be imagined~\cite{Laidler}. It has been known that the realm of subatomic particles has inherent uncertainty. This means that one can never be certain about a subatomic particle's position - or any other physical property. The subatomic world is described by probabilities. The fact that our world is not determinate at the fundamental level has been the subject of controversy not only in physics but also in philosophy. Hence, many researchers have been interested in solving this problem using various methods. One of these methods is the use of neural network methods ~\cite{schrodingernet, Schrodinger, deepschrodinger, deepschrodingepfau}. The Schrodinger equation, written by Erwin Schrodinger in 1926, is as follows ~\cite{Hall}:

\begin{equation}
    \label{shrodinger}
    \begin{split}
   i\hbar\frac{d^{2}\psi(x)}{dx^{2}}=\hat{H}\psi(x),
    \end{split}
\end{equation}

where $\psi$ is the system's wave function, $\hbar$ is Planck's constant, and $\hat{H}$ is the Hamiltonian of the system. The Hamiltonian is a mathematical operator that describes the energy of the system. All materials exhibit wave-particle duality, meaning they have properties of both waves and particles. However, this does not imply, for instance, that an electron is itself a wave. The wave functions refer to mathematical functions that reflect the likelihood of locating a particle in a specific location. In certain situations, subatomic particles act as waves. This concept was first introduced by Max Born in 1926. It explains that the probability density of finding a particle at a specific point is proportional to the square of the magnitude of the particle's wave function at that point. The wave function is complex, with time and three spatial coordinates as its parameters. It's important to note that "complex" here doesn't mean "complicated," but rather that it yields complex numbers. This complexity initially made it challenging to apply to the real world, and it's not possible to describe a particle's position or velocity with a real function.  It took some time to find an interpretation that later led Max Born to propose the idea that the wave function might be related to probabilities. This way of expressing probabilities is in terms of real numbers, and therefore, instead of directly relating its values to probabilities, we can take its squared magnitude ~\cite{Hall,deepschrodinger} - any complex number or function has a magnitude - which yields a real value and is correct!

The wave function for a particle with momentum $p$ and energy $E$ is expressed as follows:

\begin{equation}
    \label{Borglie}
    \begin{split}
   \psi(x)=e^{i(kx-w)}.
    \end{split}
\end{equation}
Subatomic particles behave differently from classical objects. They exhibit both wave-like and particle-like characteristics. In addition, the equations that describe particles at that scale have a probabilistic nature. The wave function in the Schrodinger equation is a function that gives the probability of finding a particle at a specific point in space and time. Its squared magnitude yields the probability density - measured in probability per unit volume - of finding the particle at a particular location in space. This is exactly what the wave function represents. A normalized one-dimensional wave function is as follows:

\begin{equation}
     \begin{split}
       \int_{-\infty}^{+\infty}\psi^{*}(x,t)\psi(x,t)dx=1.
    \end{split}
\end{equation}

since the squared magnitude of the wave function - which is equal to the product of the wave function and its complex conjugate inside the integral - gives the probability density, it can have a total probability of 1 or 100. If these probabilities are considered throughout the entire space, they lead to the wave function equation. Consider the Schrodinger equation \ref{shrodinger}, which takes into account everything that has been said so far, describing the probability waveforms, how they evolve, and how they behave under external influences.
The Hamiltonian operator can be represented as follows, where H is the Hamiltonian operator
\begin{equation}
\label{hamiltoni}
     \begin{split}
       \hat{H}\equiv-\frac{\hbar^{2}}{2m}\frac{\partial^{2}}{\partial x^{2}}+V(x,t).
    \end{split}
\end{equation}
According to the Schrodinger equation \ref{shrodinger} and the Hamiltonian operator, we have ~\cite{Schrodinger}

\begin{equation}
\label{timedepended}
     \begin{split}
       -i\hbar\frac{\partial}{\partial t}\psi(x,t)=(-\frac{\hbar^{2}}{2m}\frac{\partial^{2}}{\partial x^{2}}+V(x,t))\psi(x,t),
    \end{split}
\end{equation}

$V(x,t)$ is an external potential applied to the system, but the kinetic energy is not clear. The first term seems to be related to the kinetic energy of a specific state. The Hamiltonian operator depends on the total energy of the system. The Schrodinger equation expresses that the wave function - or quantum state - changes over time, and its evolution depends on the total energy - potential + kinetic - of the system.
The Schrodinger equation is one of the most famous equations in all of physics. It allows for precise predictions about various quantum systems and their time evolution. However, there are limitations to its application. It seems that even with the most powerful computers, solving the Schrodinger equation for systems with many particles is very difficult. Many experts hope that with the invention of quantum computers, this limitation will be lifted. But for now, this beautiful equation is not suitable for describing quantum systems on a large scale.
As previously mentioned, the Schrodinger equation describes how a quantum system evolves, where each particle is described by another particle.
The time-dependent Schrodinger equation (TDSE) allows for the possibility of stationary waves, we can solve for them specifically by simplifying the TDSE into the time-independent Schrodinger equation (TISE). It is assumed that the temporal and spatial parts of the solution can be obtained through variable separation, and focus can be placed on the part of the solution that only involves spatial derivatives. This part of the equation is the time-independent Schrodinger equation.

\begin{equation}
\label{hamiltoni}
     \begin{split}
       E\psi(x)=(-\frac{\hbar^{2}}{2m}\frac{\partial^{2}}{\partial x^{2}}+V(x))\psi(x),
    \end{split}
\end{equation}

It is clear, that there is no time dependence. Stationary waves remain constant over time, so the probability density of the particle remains constant. Instead of a partial derivative concerning time on the left-hand side of equation \ref{timedepended}, there is a fixed energy (E) indicating the energy of the state.

The time-independent Schrodinger equation is an eigenvalue problem. This means it can be represented in matrix form as follows ~\cite{Hall, deepschrodingepfau}

\begin{equation}
\label{matrix}
     \begin{split}
       H\psi=E\psi,
    \end{split}
\end{equation}
In this equation, $H$ is the Hamiltonian matrix (which is fundamentally the sum of the kinetic energy and the potential energy of a particle), $\psi$ is the wave function vector, and $E$ is the eigenvalue of energy. This relationship means that multiplying the matrix $H$ by the vector $\psi$ yields the same result as multiplying the scalar value $E$ by the vector $\psi$, and physically, this implies that the Hamiltonian operation, which is the sum of the kinetic and potential energy, returns the total energy of the particle.

As seen in equation \ref{matrix}, the matrix $H$ is the sum of the second-order derivative (which is the kinetic energy in terms of physics) and the potential energy $V$. For simplicity, we assume both the Planck constant $\hbar$ and the mass of the particle $m$ are equal to 1.
This article explores the solution of the time-independent Schrodinger equation in two dimensions through the use of a Hermite neural network. Additionally, the neural network was simulated using the MATLAB simulator. First, we discuss the history of solving differential equations using different numerical methods. 
The different categories of numerical methods for solving differential equations are as follows \cite{Drparand20041}: Finite Element, Finite Difference, Spectral Methods, and Meshless.
The selection of one of these methods depends on the type of equation and the specific domain in which the problems are defined. Spectral methods demonstrate excellent performance for problems characterized by a relatively smooth and regular geometry. They are highly efficient and accurate, encompassing methods like collocation, Galerkin, Petro-Galerkin, and Tau \cite{Boyd2000, Drparand20041,shen201115}. It uses the orthogonal basis that is the solutions of the Sturm-Liouville equation. The basis exhibits specific behaviors depending on the type of polynomial they have. The choice of basis function is an important feature of spectral methods, which focuses on it. Some problems are in semi-infinite or infinite intervals, so for numerically solving problems in these intervals, orthogonal polynomials such as Hermite, Laguerre, and sine functions can be used. There are also other methods for solving these types of problems such as mapping Chebyshev, and Legendre functions in infinite intervals, as well as cutting semi-infinite intervals or transforming the problem in a semi-infinite interval to a finite interval problem using a variable change. As mentioned, spectral methods can be considered an extension of weighted residual methods. Based on this, we will briefly explain weighted residual methods and choices of weight function in the next section. There is another method for solving differential equations, which has attracted the attention of researchers in the last few years – the use of machine learning algorithms \cite{ordinaryneuranetwork, Guo}. Machine learning algorithms essentially function as approximation functions. When trained on a dataset of inputs and outputs, they compute a mathematical function or a set of operations to relate the inputs to the outputs. 
Machine learning algorithms are typically trained to approximate functions mapping inputs to Euclidean space outputs using classical graphs with x, y, and z axes. However, a new approach defines inputs and outputs in Fourier space. Since Fourier approximation is much easier in Fourier space compared to solving differential equations in Euclidean space, working with machine learning algorithms becomes more convenient. As mentioned above, most differential equations exist in semi-infinite or infinite intervals. The significance of this fact has led to the realization.  To improve numerical methods for solving differential equations, they have used different basis functions \cite{drparand2019,drparand2018-125,drparand2019/111}

Parand \cite{drparand2022-2, drparand2021} solved different differential equations by using supervised and unsupervised machine learning methods and combining them with spectral methods. Parand \cite{Drparand20215} was able to provide a new solution to optimize this algorithm by using spectral methods and combining them with neural networks.

\section{Implement}

Initially, we will discuss the method for solving this problem. Next, we will provide an overview of the neural network structure used in this article. Finally, we will present the proposed solution to the Schrodinger equation.
\subsection{Hermite functions}
In this section, we consider the properties of Hermite functions.
$\widetilde{H}_{n}(x)$ are the normalized Hermite functions of degree $n$, which describe the properties of Hermite functions \cite{shen,drParand2018-2}.
\begin{equation}
\widetilde{H}_{n}(x)=\frac{1}{\sqrt{2^{n}n!}}e^{\frac{-x^{2}}{2}}H_{n}(x),\quad n\geq0,\:x\in\mathbb{R}.
\end{equation}

The formula for orthogonal relation for Hermite functions is as follows:

\begin{equation}
\int_{-\infty}^{+\infty}\widetilde{H}_{n}(x)\widetilde{H}_{m}(x)=\sqrt{\pi}\delta_{mn},
\end{equation}

where $\delta_{mn}$is the Kronecker delta function. Hermite functions have a recurrent relation defined in the  $(-\infty,+\infty)$ domain

\begin{gather}
\widetilde{H}_{n+1}(x)=x\sqrt{\frac{2}{n+1}}\widetilde{H}_{n}(x)-\sqrt{\frac{n}{n+1}}\widetilde{H}_{n-1}(x),\quad n\geq1,\nonumber \\
\widetilde{H}_{0}(x)=e^{\frac{-x^{2}}{2}},\;\widetilde{H}_{1}(x)=\sqrt{2}xe^{\frac{-x^{2}}{2}}.
\end{gather}

Use the Hermite functions' recurrence relation and formula to get the result

\begin{equation}
\widetilde{H^{\prime}}_{n}(x)=\sqrt{2n}\widetilde{H}_{n-1}(x)-x\widetilde{H}_{n}(x)=\sqrt{\frac{n}{2}}\widetilde{H}_{n-1}(x)-\sqrt{\frac{n+1}{2}}\widetilde{H}_{n+1}(x),
\end{equation}

and it becomes

\begin{equation}
\int_{-\infty}^{+\infty}\widetilde{H^{\prime}}_{n}(x)\widetilde{H^{\prime}}_{m}(x)dx=\begin{cases}
-\frac{\sqrt{n\pi(n-1)}}{2}, & m=n-2,\\
(n+\frac{1}{2})\sqrt{\pi}, & m=n,\\
-\frac{\sqrt{\pi(n+1)(n+2)}}{2}, & m=n+2,\\
0, & Otherwise.
\end{cases}
\end{equation}

\begin{equation}
\tilde{P}:\{u:u=e^{\frac{-x^{2}}{2}}\nu,\;\forall_{\nu}\epsilon P_{N}\},
\end{equation}
    
 Where $P_{N}$ represents the Hermite polynomials of degree $N$.   
    
\subsection{Solution method}
Consider the ordinary differential equation of the following form
\begin{equation}
\label{diffnn}
     \begin{split}
\begin{array}{c}
f(x,g(x),g^{'}(x),g^{''}(x,\ldots,g^{n}(x))=0, \quad x\in R,
\end{array}
    \end{split}
    \end{equation}

    where $g(x)$ is the function to be found and $g^{n}(x)$ is the nth derivative of the function $g(x)$. The trial solution of $g(x)$ will be as follows

     \begin{equation}
     \begin{split}
\begin{array}{c}
g_{t}(x)=h_{1}(x)+h_{2}(x,N(x,P)),
\end{array}
    \end{split}
    \end{equation}

$h_{1}(x)$ is a function that satisfies a set of conditions $g_{t}(x)$ and $N(x,P)$ is a neural network described by $P$ with weights and biases. The role of $h_{2}(x, N(x, P)$) is to show that the output of $N(x, P)$ is zero when $g_{t}(x)$ satisfies the conditions for the values of $x$.
The reason why the neural network is used, as explained earlier, is an optimization method to minimize parameters, weights, and bias, and in the proposed method, it is done through backward propagation. To define the minimization, a cost function must be defined to minimize it, so we set the equation \ref{diffnn} equal to zero. We can consider the mean squared error as a cost function for the input $x$. The cost function $c(x),P)$ is defined as follows

\begin{equation}
     \begin{split}
\begin{array}{c}
c(x,P)=(f(x,g(x),g^{'}(x),g^{''}(x),\ldots,g^{n}(x))^2.
\end{array}
    \end{split}
    \end{equation}
    
If $N$ input is given as vector $x$ with elements $x_{i}\;i=1,\ldots,N$, the cost function will be defined as follows

    \begin{equation}
    \label{cost}
     \begin{split}
\begin{array}{c}
c(x,P)=\frac{1}{N}\sum_{i=1}^{N}(f(x,g(x),g^{'}(x),g^{''}(x),\ldots,g^{n}(x))^{2}.
\end{array}
    \end{split}
    \end{equation}

In the neural network, the parameter $P$ is set in such a way that the cost function of the equation \ref{cost} is minimized. In this present, minimization is done by the gradient descent in the equation \ref{cost} and there are various libraries in Python for numerical derivatives. Autograd is used in this paper.
The next step is changing the parameters to minimize the cost function.
If $\overrightarrow{x}$ is considered as a vector with elements $x_{i}, i=1, \ldots,N$, the absolute or squared difference should be near zero, ideally zero\cite{introductionneuralnetwork,neuralnetwork}.

 \begin{equation}
     \begin{split}
\begin{array}{c}
\begin{array}{c}
c(\overrightarrow{x},P)=\frac{1}{N}\sum_{i}(g_{t}^{'}((x_{i}),P)-F(N((x),P))^{2}),\\
g_{t}^{'}(x),P)=A(x)+F(N((x),P)),
\end{array}
\end{array}
    \end{split}
    \end{equation}
    
    To minimize the cost function, an optimization method should be selected. As mentioned before, in this paper, gradient descent is used.
     The idea of the gradient descent algorithm is to update the parameters in the direction that the cost function is minimized. This method finds the optimal boundary during different iterations. In this way, it begins from one place and moves in the direction of the negative slope of the error, and when the slope of the error becomes zero (minimum error), the training process stops and does not continue.
     In general, updating some parameters $\overrightarrow{w}$ according to a defined cost function $c(\overrightarrow{x},\overrightarrow{w})$ is as follows

            \begin{equation}
     \begin{split}
\begin{array}{c}
\mathbf{w}_{new}=\mathbf{w}-\lambda\nabla_{w}c(\mathbf{x},\mathbf{w}),
\end{array}
    \end{split}
    \end{equation}
The number of iterations takes place until it is smaller than $ \big|\big| \boldsymbol{\omega}_{\text{new} } - \boldsymbol{\omega} be \big|\big|$. The value of $\lambda$ determines the algorithm's steps in the direction $ \nabla_{\boldsymbol{\omega}} C(\boldsymbol{x}, \boldsymbol{\omega})$.
$\nabla_{w}$ gradient sign is expressed according to the elements in $w$ and the cost function $c( x, P)$ is calculated using the chain derivative, which causes the values of the weights to be obtained and ultimately the error values to be minimized.

    Next, we will discuss the method of solving partial differential equations, which involves the same techniques as solving ordinary differential equations. Generally, if we have a function $g(x_{1},\ldots,x_{N})$ with N variables, its partial differential equation is expressed as follows:

    \begin{equation}
     \begin{split}
f(x_{1},\ldots,x_{N},\frac{\partial g(x_{1},\ldots,x_{N})}{\partial x_{1}},\ldots,\frac{\partial g(x_{1},\ldots,x_{N})}{\partial x_{N}},\frac{\partial g(x_{1},\ldots,x_{N})}{\partial x_{1}\partial x_{2}},\ldots,\frac{\partial^{n}g(x_{1},\ldots,x_{N})}{\partial x_{N}^{n}}),
    \end{split}
    \end{equation}
    where $f$ includes derivatives from $g(x_{1},\ldots,x_{N})$ up to order $n$. The trial solution can be expressed in the following form

\begin{equation}
     \begin{split}
     \label{neuralpartical}
g_{t}(x_{1},\ldots,x_{N})=h_{1}(x_{1},\ldots,x_{N})+h_{2}(x_{1},\ldots,x_{N},N(x_{1},\ldots,x_{N},P)),
    \end{split}
    \end{equation}
    $h_{1}(x_{1},\ldots,x_{N})$ satisfies some conditions for $g_{t}(x_{1},\ldots,x_{N})$ . The neural network $N(x_{1},\ldots,x_{N},P)$ has weights and biases that are described by $P$ and you get the output of the network $h_{2}(x_{1},\ldots,x_{N},N(x_{1},\ldots,x_{N},P))$.
     As previously discussed the cost function, the cost function is the average squared error that the network should try to minimize. Here, the aim is to minimize the cost function. To minimize equation \ref{neuralparticle}, $P$ needs adjustment, considering several variables.

           \begin{equation}
     \begin{split}
\begin{array}{c}
c(x_{1},\ldots,x_{N},P)=(f(x_{1},\ldots,x_{N},\frac{\partial g(x_{1},\ldots,x_{N})}{\partial x_{1}},\ldots\\
\ldots,\frac{\partial g(x_{1},\ldots,x_{N})}{\partial x_{N}},\frac{\partial g(x_{1},\ldots,x_{N})}{\partial x_{1}\partial x_{2}},\ldots,\frac{\partial^{n}g(x_{1},\ldots,x_{N})}{\partial x_{N}^{n}}))^{2},
\end{array}
    \end{split}
    \end{equation}
when considering a set of values for x, represented as $\mathbf{x}=(x_{1},\ldots,x_{N})$, the expression for the cost function can be presented as follows
    \begin{equation}
     \begin{split}
\begin{array}{c}
c(\mathbf{x},P)=(f(\mathbf{x},\frac{\partial g(x_{1},\ldots,x_{N})}{\partial x_{1}},\ldots,\frac{\partial g(x_{1},\ldots,x_{N})}{\partial x_{N}},\\
\ldots,\frac{\partial g(x_{1},\ldots,x_{N})}{\partial x_{1}\partial x_{2}},\ldots,\frac{\partial^{n}g(x_{1},\ldots,x_{N})}{\partial x_{N}^{n}}))^{2},
\end{array}
    \end{split}
    \end{equation}
if we have different sets  $M$  and values  $x_{1},\ldots,x_{N}$, where $x_{i}=(x_{1}^{(i)},\ldots,x_{N}^{(i)}),\;i=1,\ldots,M$ represents the rows of the matrix $X$, then the cost function can be reformulated to another shape.

   \begin{equation}
     \begin{split}
c(\mathbf{X},P)=\sum_{i=1}^{M}(f(\mathbf{x_{i}},\frac{\partial g(x_{i})}{\partial x_{1}},\ldots,\frac{\partial g(x_{i})}{\partial x_{N}},\ldots,\frac{\partial g(x_{i})}{\partial x_{1}\partial x_{2}},\ldots,\frac{\partial^{n}g(x_{i})}{\partial x_{N}^{n}}))^{2},\ldots \;i=1,\ldots,M.
    \end{split}
    \end{equation}
At present, it has been tried to use the Adam algorithm because one of the problems of some descending algorithms is that they face local minima and may get caught in this trap. It is an optimization algorithm that can be used instead of the classical stochastic gradient descent method to update the weights of the iterative network based on the training data. This method was presented by Diederik Kingma in 2015. Its name is derived from Adaptive Moment Estimation. In stochastic gradient descent, there is a learning rate that does not change during learning, but in this method, there is a learning rate for each weight of the network, and it is adapted separately with the expansion of learning.

\subsection{Hermite Neural network structure}

\begin{figure}[h]

\begin{centering}
\includegraphics[scale=0.5]{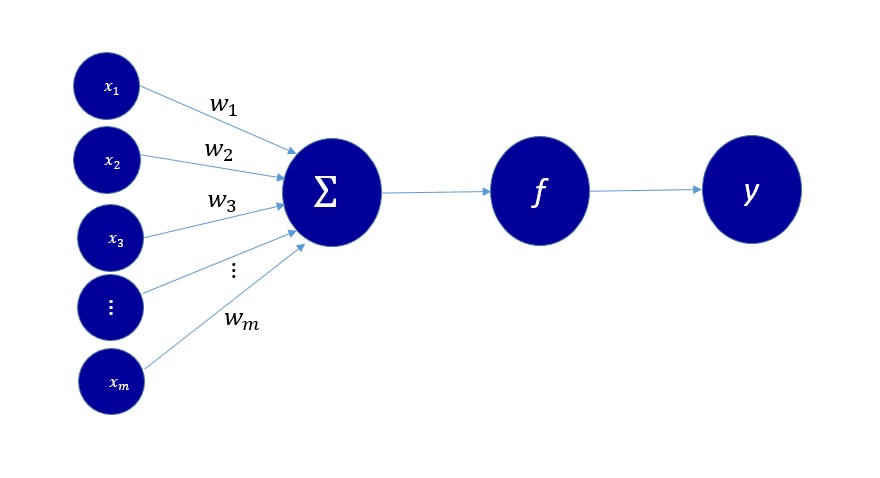}
\par\end{centering}
\caption{Neural Network Structure, $x_{i}:(i=1,...,m)$ are the input vectors, $w_{i}:(i=1,...,m)$ are the weight vectors, $\sum$ is the sum operator, $f$ is the activation function and $y$ is output. }\label{neural}
\end{figure}

The mathematical neuron model is a simulation of the biological neuron. A nerve cell contains an array of dendrites that receive signals from the environment. Dendrites, in other words, are our inputs. The inputs are represented in the mathematical model by the vector $X_{m\times1}$, where $m$ is the number of inputs, as illustrated in the figure below. The chemical interactions that occur for the inputs are represented by the coefficients assigned to these inputs, which, as shown in the figure \ref{neural}, are expressed as $w_i:(i=1,\dots,m)$ and are referred to as weights. These chemical interactions, for example, might decrease or amplify the signal, which can be approximated with a factor smaller or greater than one. Furthermore, these chemical interactions might result in the figure \ref{neural}.

As shown in the image above, the result of $\sum_{i=1}^{m}w_{i}x_{i}$ is placed into the cell's core, which can be named $Z$, and the core determines whether or not this value is greater than a threshold limit (Bias) represented by $b$. If $Z$ is larger than $b$, its output is activated; otherwise, the number zero or the negative number one is displayed in the output, indicating that it is inactive, as illustrated in the figure by the function $f$. As a result, its ultimate output is $y=f(Z)$. The Activation Function is the name given to this function. According to the above explanations, we have the following relations

 \begin{equation}
     \begin{split}
if\quad Z=\sum_{i=1}^{m}w_{i}x_{i}\geq b \implies y=f(z)=1,\\if\quad Z=\sum_{i=1}^{m}w_{i}x_{i}<b \implies y=f(z)=-1,
    \end{split}
    \end{equation}
    
    that the result of $\sum_{i=1}^{m}w_{i}x_{i}$ can also be displayed by multiplying two vectors $X$ and $W$
    \begin{equation}
     \begin{split}
\mathbf{X}=\left[\begin{array}{c}
x_{1}\\
x_{2}\\
\vdots\\
x_{m}
\end{array}\right]=\left[\begin{array}{cccc}
x_{1} & x_{2} & \ldots & x_{m}\end{array}\right]^{T}\qquad\mathbf{W}=\left[\begin{array}{c}
w_{1}\\
w_{2}\\
\vdots\\
w_{m}
\end{array}\right]=\left[\begin{array}{cccc}
w_{1} & w_{2} & \ldots & w_{m}\end{array}\right]^{T},
    \end{split}
    \end{equation}

   Hence, it is stated

     \begin{equation}
     \begin{split}
\theta^{T}X_{a}=\left[\begin{array}{ccccc}
w_{1} & w_{2} & \ldots & w_{m} & b^{'}\end{array}\right]^{T}\left[\begin{array}{c}
x_{1}\\
x_{2}\\
\vdots\\
x_{m}\\
1
\end{array}\right]=\begin{array}{ccccc}
w_{1}x_{1}+ & w_{2}x_{2}+ & \cdots+ & w_{m}x_{m}+ & b^{'}=\mathbf{W}^{T}\mathbf{X}+b^{'},\end{array}
    \end{split}
    \end{equation}

  $\theta^{T}X_{a}$ is $\mathbf{W}^{T}\mathbf{X}+b^{'}$ or $\mathbf{W}^{T}\mathbf{X}-b$. Inequalities related to the function $f(Z)$, can be displayed in the following form

\begin{equation}
     \begin{split}
if\quad\theta^{T}X_{a}\geq0\implies y=1,\\if\quad\theta^{T}X_{a}<0\implies y=-1,
    \end{split}
    \end{equation}
    
When we consider the vector $X$ as an input consider the calculations related to the weight functions and apply the active function for each neuron, This process is carried out from the middle layer to the output layer, which is called Feedforward neural networks. But if we also consider reverse connections, i.e. we want the connection of one layer with our previous layer, in this case, it will require a more complex operation. One of the methods used for optimization is the use of gradient descent. In the figure \ref{neural}, we can see the structure of a multi-layer neural network. The input layer with $m$ neurons is $X=[x_{1},x_{2},...x_{m}]\subset\mathbb{R}^{m}$ and the output layer is $Y\subset\mathbb{R}^{n}$ \cite{drparand2022a}.

  \begin{figure}[h!]

\begin{centering}
\includegraphics[scale=0.4]{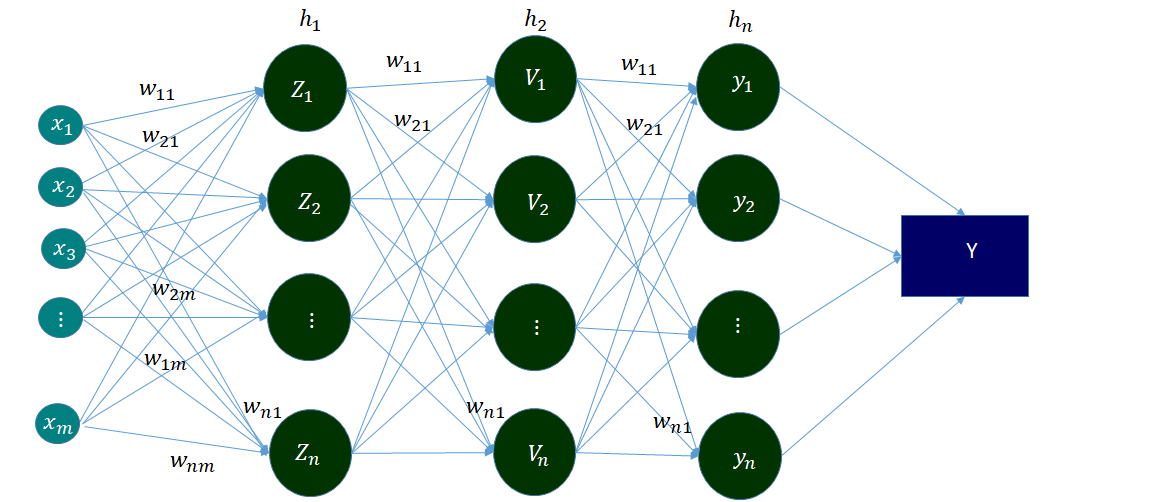}
\par\end{centering}
\caption{Feedforward}\label{neurall}

\end{figure}

An artificial neural network according to the form \ref{neurall} can be represented as follows

 \begin{equation}   
\left[\begin{array}{c}
z_{1}\\
z_{2}\\
z_{3}\\
\vdots\\
z_{n}
\end{array}\right]=\left[\begin{array}{cccc}
w_{11}^{(1)} & w_{21}^{(1)} & \cdots & w_{n1}^{(1)}\\
w_{12}^{(1)} & w_{22}^{(1)} & \cdots & w_{n2}^{(1)}\\
w_{13}^{(1)} & w_{23}^{(1)} & \cdots & w_{n3}^{(1)}\\
\vdots & \vdots & \vdots & \vdots\\
w_{1m}^{(1)} & w_{2m}^{(1)} & \cdots & w_{nm}^{(1)}
\end{array}\right]\left[\begin{array}{c}
x_{1}\\
x_{2}\\
\vdots\\
x_{m}
\end{array}\right]\text{,}
\end{equation}

where $Z=[z_{1},z_{2},...z_{n}]^{T}\subset\mathbb{R}^{n}$ is the output vector resulting from the dot multiplication of the hidden layer weights vector $h_1$ and the input vector $X=[x_{1},x_{2},...x_{m}]\subset\mathbb{R}^{m}$ and it can be shown in the following form

\begin{equation}
Z_{j}=W_{ij}.X_{i},
\end{equation}
$W_{ij}$ is the weights matrix, $X_{i}$ is the input vector and $Z_{j}$ is the output vector. This relationship can be written for other layers. A point that should be noted is that from this layer onward, better results can be obtained by applying different active functions.

\begin{equation}
\left[\begin{array}{c}
v_{1}\\
v_{2}\\
v_{3}\\
\vdots\\
v_{k}
\end{array}\right]=\left[\begin{array}{cccc}
w_{11}^{(2)} & w_{21}^{(2)} & \cdots & w_{n1}^{(n)}\\
w_{12}^{(2)} & w_{22}^{(2)} & \cdots & w_{n2}^{(2)}\\
w_{13}^{(2)} & w_{23}^{(2)} & \cdots & w_{n3}^{(2)}\\
\vdots & \vdots & \vdots & \vdots\\
w_{1k}^{(2)} & w_{2k}^{(2)} & \cdots & w_{nk}^{(2)}
\end{array}\right]\left[\begin{array}{c}
z_{1}\\
z_{2}\\
\vdots\\
z_{n}
\end{array}\right]\text{.}
\end{equation}
The output vector of the chosen hidden layer is created by combining the above equations, using the dot multiplication of the combined matrix and the input vector, and continuing until the output layer.

\begin{equation}
V_{k}=W_{jk}.Z_{j}=W_{ij}.W_{jk}.X_{i}=W_{ijk}.X_{i},
\end{equation}

\begin{equation}
\left[\begin{array}{c}
y_{1}\\
y_{2}\\
y_{3}\\
\vdots\\
y_{l}
\end{array}\right]=\left[\begin{array}{cccc}
w_{11}^{(n)} & w_{21}^{(n)} & \cdots & w_{n1}^{(n)}\\
w_{12}^{(n)} & w_{22}^{(n)} & \cdots & w_{n2}^{(n)}\\
w_{13}^{(n)} & w_{23}^{(n)} & \cdots & w_{n3}^{(n)}\\
\vdots & \vdots & \vdots & \vdots\\
w_{1l}^{(n)} & w_{2l}^{(n)} & \cdots & w_{nl}^{(n)}
\end{array}\right]\left[\begin{array}{c}
v_{1}\\
v_{2}\\
\vdots\\
v_{n}
\end{array}\right]\text{.}
\end{equation}

Finally, the output matrix is determined as follows

\begin{equation}
Y_{l}=W_{kl}.V_{k}=W_{ijk}.W_{kl}.X_{i}=W_{ijkl}.X_{i},
\end{equation}

From this stage, network training is backward. In short, it can be said that the training of the network is forward in one stage, that is, it goes from the inputs to the output, and in the next stage, it is backward, in which it goes backward using the computed error. Optimizing methods can be used to get better results
     In real problems, the number of data is large, the above method cannot be used.
     In this case, better results can be achieved by applying appropriate activation functions and reducing the cost function error using optimization methods. Since derivation is needed in optimization methods, derivable functions must be used. As explained in the previous section, the selection of basis functions is one of the characteristics of spectral methods, and the basis functions in spectral methods are infinitely differentiable general functions\cite{ordinaryneuranetwork,neuralnetwork}. Parand solved several orthogonal expressions \cite{drparand2022-2,drparand2021fred,drparand2021fcdn} in neural networks to solve differential equations, whose features include fast convergence and easy calculation. In this treatise, it has been tried to use Hermite functions and fractional Hermite functions in neural networks to solve differential equations.

For each continuous function $y:\left[a,b\right]\rightarrow\mathbb{R}$ and $w_{n}\;(n=0,1,2,\ldots,N)$ are the weights and $N$ is a natural value. So the neural network of Hermite functions with $N+1$ neurons and $\tilde{H}_{n}(x)$ of Hermite functions will be defined as follows

    \begin{equation}
     \begin{split}
yLNN(x)=\sum_{n=0}^{N}w_{n}\tilde{H}_{n}(x),
    \end{split}
    \end{equation}
    $yLNN$ is an approximation of $y$

        \begin{equation}
     \begin{split}
\parallel y(x)-yLNN(x)\parallel=\parallel y(x)-\sum_{n=0}^{N}w_{n}\tilde{H}_{n}(x)\parallel<\varepsilon.
    \end{split}
    \end{equation}

    Neural networks based on Hermite functions or fractional Hermite functions include three layers: input layer, hidden layer based on Hermite functions or fractional Hermite functions, and output layer. The output of this neural network model for differential equations is expressed as follows

        \begin{equation}
     \begin{split}
yLNN(x)=\sum_{n=0}^{N}w_{n}\tilde{H}_{n}(x).
    \end{split}
    \end{equation}
    Differential equations and boundary conditions can be defined as follows

        \begin{equation}
     \begin{split}
\begin{array}{c}
\mathcal{L}(y(x))=f(x),\quad x\in\Omega\subseteq\mathbb{R},\\
\beta(y(x))=\alpha,
\end{array}
    \end{split}
    \end{equation}
    $\mathcal{L}$, $\beta$ are differential equation operators, $y(x)$ is the unknown value, $d$ is the dimension of the feature space, and $f(x)$ is linear or non-linear. $\alpha$ indicates the boundary value in the specified interval, which can be a constant value. By placing the approximate solution in the formula and applying the boundary conditions, the values of the weights $W$ can be obtained, and hence the new equation will be as follows

        \begin{equation}
     \begin{split}
\begin{array}{c}
\mathcal{L}(yLNN(x))=f(x),\quad x\in\Omega\subseteq\mathbb{R},\\
\beta(yLNN(x))=\alpha,
\end{array}
    \end{split}
    \end{equation}
    By placing Hermite collocation points $z_{i},i=0,\ldots,M$ which are considered as the roots of $H_{n+1}(z)$, the above equation will be written as follows

        \begin{equation}
     \begin{split}
\begin{array}{c}
\left[\begin{array}{c}
\mathcal{L}(\sum_{n=0}^{N}\tilde{H}_{n}(z_{i}))\\
\ldots\\
\beta(\sum_{n=0}^{N}\tilde{H}_{n}(z_{boundary}))
\end{array}\right]\left[W\right]=\left[\begin{array}{c}
f_{i}\\
\ldots\\
\alpha
\end{array}\right]\quad i=0,\ldots,m,
\end{array}
    \end{split}
    \end{equation}
    where $f_{i}=f(z_{i})$ and $W=[w_{0},\ldots w_{m}]$ will be the boundary conditions $[z_{0},\ldots,z_{m}]$. This method is briefly shown in the figure 2. $X$ is the input layers and $Z=[z_{0},\ldots,z_{m}]$ is the collocation points. Functional development block is based on Hermit functions

\begin{figure}[h]

\begin{centering}
\includegraphics[scale=0.3]{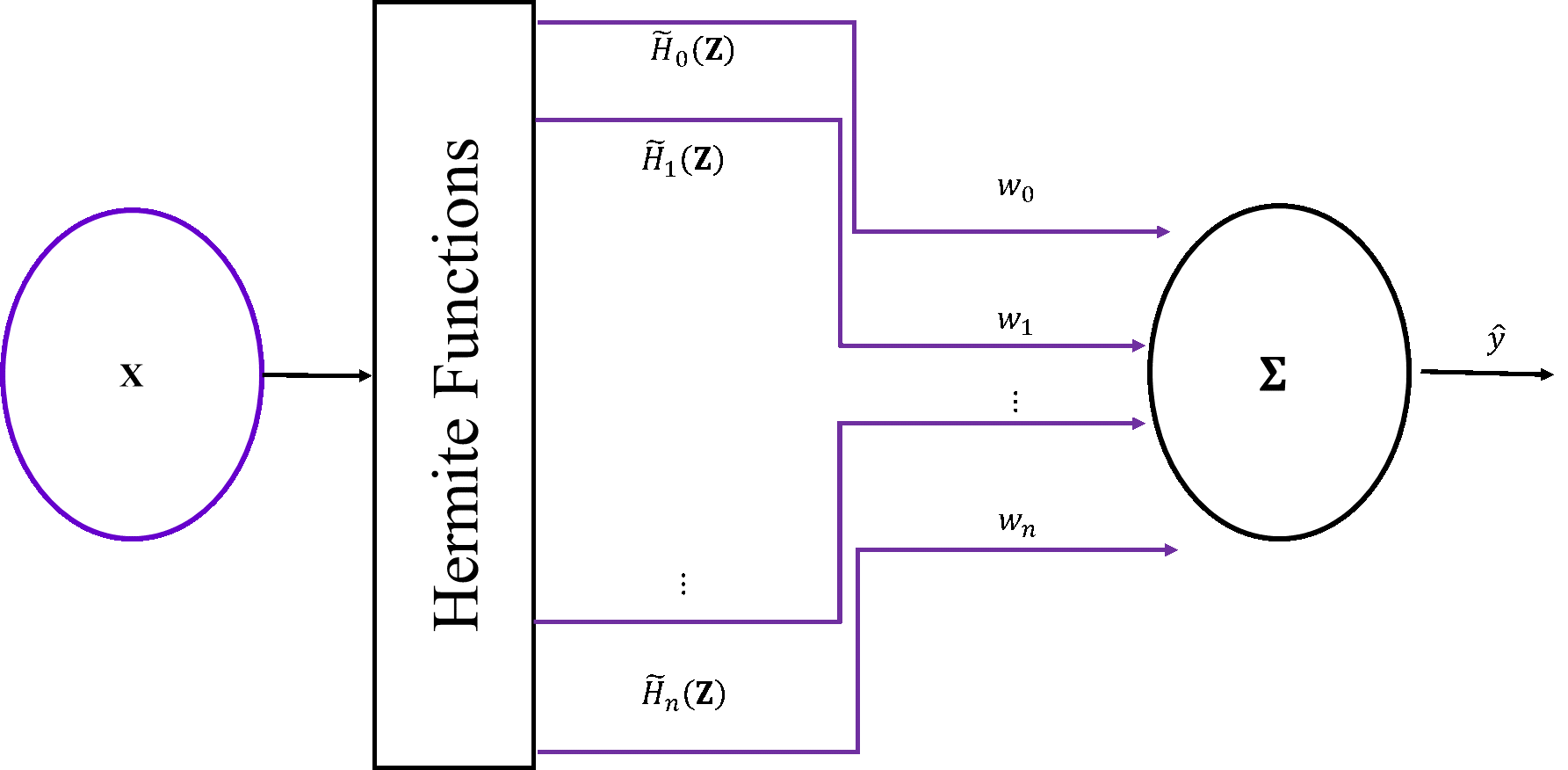}
\par\end{centering}
\begin{centering}
\caption{Hermite Neural Network Structure }\label{2}
\par\end{centering}
\end{figure}

\subsection{Solving the Schrodinger equation}

This paper defines and trains a neural network to solve the Schrodinger equation in an infinite well potential. First, various constants are defined, including the dimensions of the computational region, the neural networks, and the weights and validity for each layer.
Then the activation function for the network is supplied (which are the Hermite functions here). The neural network is trained using forward learning. In each step, several matrices and vectors containing the available approximate values are prepared, and then the error is calculated and the various weights and credits are updated.
Finally, various graphs are obtained to show the results, including energy, wave function, and neural network graphs.
In the following section, we will describe the method implementation and show the simulation structure for this method. Using this simulator, we could show the results of solving this equation. 
One of the important parts of this method is the determination of training points, where the roots of Hermite functions are used for training.  Using these roots, we can determine the points at which the positive and negative repetitions of the Hermitian functions occur. These points are used as training points and the neural network inputs are placed in these points. Through this method, the neural network can be trained on the Hermite roots and continue improving its model between them. By using this approach, the neural network can approach points near the Hermite roots and thus have more training accuracy and efficiency.

Input layer: This layer has two nodes, each of which is equal to $x$ and $y$ coordinates. These two nodes represent the input data to the network.

Hidden layer: This layer has 10 hidden layers and there are 5 nodes in each layer. The Hermite activation function is used to calculate the output of each node of this layer.

Output layer: This layer has only one node that represents the output of the network. In this code, this node is the desired output for estimating the wave mode.

In network training, using the error back-propagation algorithm, the weights and biases are updated in each step and the network tries to estimate the desired wave mode more accurately.

The domain we have defined in $X,Y$ is in a specific interval that we have used based on the definition of the problem. The activation function used the Hermite functions.   

The steps for computing approximate values and updating weights during training are as follows:

The output of the first hidden layer
\begin{equation}
z_{1}=W_{1}^{'}*[\psi_{x}^{'};\psi_{y}^{'}]+b_{1},
\end{equation}
where  $W_{1}$ is the weight matrix between the input layer and the first hidden layer, $\psi(x,y)$  is the wave function for the input variables $x$ and $y$, and $b_{1}$ is the bias vector of the first hidden layer. We define the activation function to obtain the output of the first hidden layer as follows:
\begin{equation}
a_{1}=yLNN(z_{1}),
\end{equation}
$yLNN$ is the approximate solution based on the Hermite activation function. The output of the second hidden layer

\begin{equation}
z_{2}=W_{2}^{'}*a_{1}+b_{2},
\end{equation}
where $W_{2}$  is the weight matrix between the first hidden layer and the second hidden layer, $a_{1}$ is the output of the first hidden layer,  and $b_{2}$ is the bias vector of the first hidden layer.  Through the internal multiplication between the weight matrix and the output of the first hidden layer $(a_{1})$, by adding the bias vector and applying the activation function on this output $z_{2}$ we get the final output of the second hidden layer $a_{2}$. For the second hidden layer, we use the Hermite activation function and define it as follows
\begin{equation}
a_{2}=yLNN(z_{2}).
\end{equation}
The computation of the output of the last layer $z_{3}$, which is an approximate value for the wave function, will also be as follows:

\begin{equation}
z_{3}=W_{3}^{'}*a_{2}+b_{3},
\end{equation}

where $W_{3}$ is the weight matrix between the second hidden layer and the output layer (final layer), $a_{2}$ is the output of the second hidden layer and $b_{3}$ is the bias vector of the output layer. By using the internal multiplication between the weight matrix $W_{3}$ and the output of the second hidden layer $a_{2}$, then by adding the bias vector $b_{3}$ and applying the activation function on it, we get the final output of the neural network in the form of vector $z_{3}$. The final output for the wave function will be

\begin{equation}
\bar{\psi}=z_{3}.
\end{equation}
We have defined the error using the difference between the actual solution $\psi$ and the predicted solution $\bar{\psi}$. First, we double the error value and multiply it in the previous layers. The reason is that multiplying by 2 can help us obtain appropriate changes in weights.

In the back-propagation algorithm, an error function is called the "loss function". This function is used to measure the amount of error of the predicted with the actual values.

We compute the partial derivative of the error function for the weight and bias variables in the network using the chain rule, so it will be
\begin{equation}
\delta_{3}=2*(\psi-\bar{\psi),}
\end{equation}
where $\delta_{3}$ is the partial derivative of the error to the weight variables the actual $\psi$ and the predicted $\bar{\psi}$. The estimate of the error in the second hidden layer $\delta_{2}$ is obtained by multiplying the error of the previous layer by the weights of the second hidden layer and also the bias of the second hidden layer.

\begin{equation}
\delta_{2}=W_{3}*(a_{2}*(1-a_{2}))*\delta_{3}.
\end{equation}

computation of the error in the first hidden layer $\delta_{1}$ by multiplying the error of the previous layer by the weights of the first hidden layer and also the bias of the first hidden layer

\begin{equation}
\delta_{1}=W_{2}*(a_{1}-(1-a_{1}))*\delta_{2}.
\end{equation}

In this method, the Stochastic Gradient Descent is used, and the gradients computed based on random samples of the data are updated with weights and biases. That is, instead of using all the data to calculate the gradient and update, a random sample of the data is used. Due to the limited use of data, this method is faster in training time and can be useful in large data sets. Then, using the mean square error $(MSE)$ between the predicate $\bar\psi$ and the actual solution $\psi$, the gradients are calculated using the chain rule. Finally, the weights and biases are updated using the learning rate $(LR)$ and gradients.

 The input data set consists of the x and y positions computed for each point, the wave function, and the wave function approximate. Then we estimate the error by the mean square error for each point. This error shows the amount of difference between the wave function and the one estimated by the Hermite neural network.

\begin{equation}
MSE=\frac{\sum_{i}(\varphi_{i}-\bar{\varphi_{i}})^{2}}{N},
\end{equation}
where $\psi$ is the actual solution and $\bar\psi$  is the approximate solution. $N$ is The number of points in the neural network, This value shows the neural network model can estimate the wave function correctly.

\begin{figure}[h]

\begin{centering}
\includegraphics[scale=0.5]{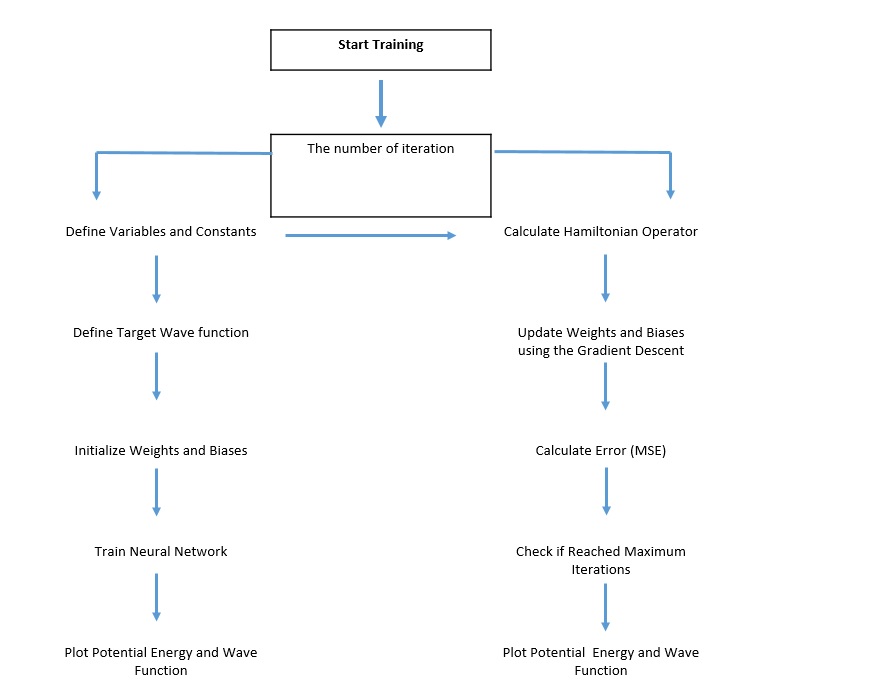}
\par\end{centering}
\begin{centering}
\caption{The present method flowchart }\label{4}
\par\end{centering}
\end{figure}

You can solve differential equations graphically and more easily using Simulink in MATLAB. This approach is very useful because it is possible to interact with the model and change the parameters in Simulink and you can easily analyze the results. The flowchart below describes the steps below in the figure \ref{4}

\begin{figure}[h]

\begin{centering}
\includegraphics[scale=0.5]{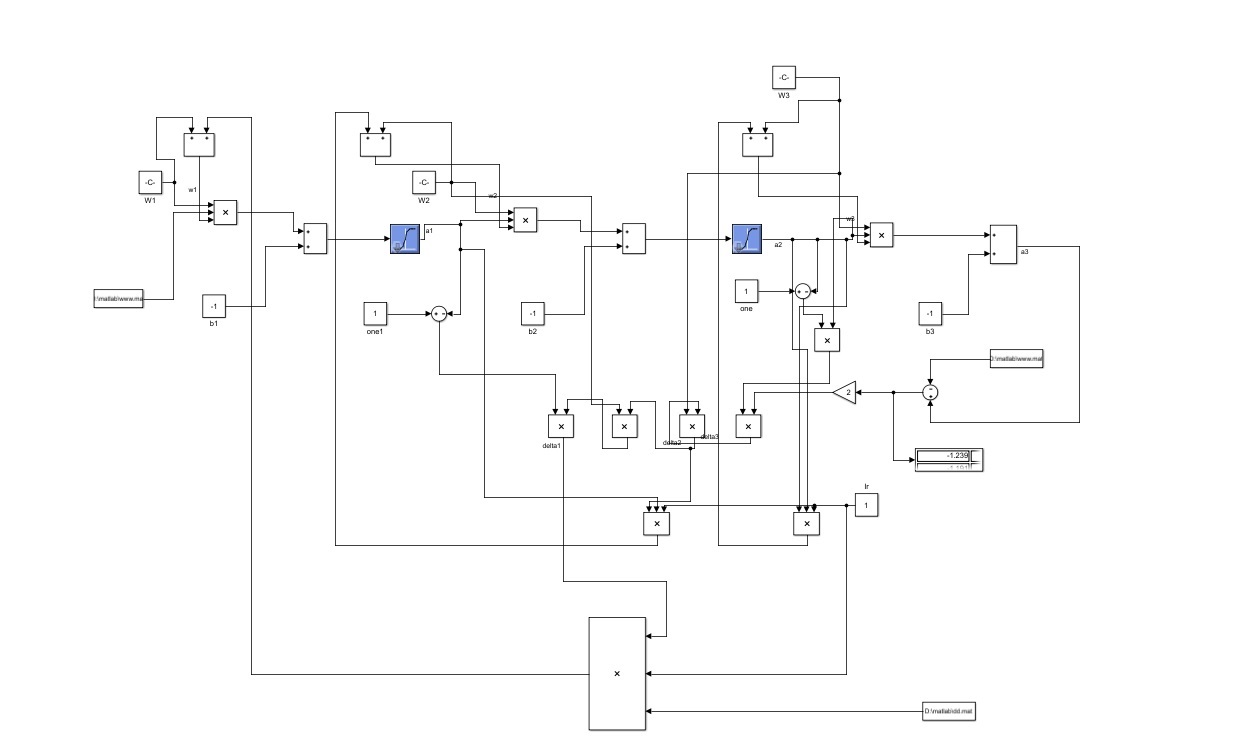}
\par\end{centering}
\begin{centering}
\caption{Neural Network Structure }\label{1}
\par\end{centering}
\end{figure}

In the next section, we consider Simulink this method. The Neural network structure simulated in the figure
\ref{1}

Input Layer: In the Simulink model, add a  From Workspace block and provide our inputs to this block. This block defines the input layer of the network.

Hidden layers: For this part, choose from the Transfer Function or Sum block and specify the number of neurons and the corresponding activation functions.

Determine and connect weights (Weights): We used the Gain block to connect layers. This block defines weights between layers.

Determining and Connecting Biases: The Sum block connects biases between layers. This block defines biases between layers.

Activation Functions: Set the activation functions of the layers in the corresponding blocks. In this example, we used the Transfer Function block for showing Hermite functions.

   \section{Numerical Results} 

The Schrodinger equation is a fundamental concept in quantum physics, used to describe the different states of a quantum system. It can be represented by the following equation \cite{Atkins}:

\begin{equation} \hat{H} = -\frac{{\hbar^2}}{{2m}} \nabla^2 \psi + V \psi + \frac{1}{2} m \omega^2 ( X^2 + Y^2 ) \psi, \end{equation}

The equation $\omega$ represents the angular frequency, and $H$ is the Hamilton operator, which consists of three components. The first term represents the particle's kinetic energy, which depends on Planck's quantum constant, denoted as $\hbar$, and the mass of the particle, denoted as $m$. Operator $\nabla^2$ is the Laplacian symbol.

The second term, denoted as $V$, represents the potential between particles. This potential can vary with space, time, and other factors such as internal and external potentials.

The third term describes the potential energy of a particle in a harmonic field or potential, modified by the function $\psi$.

\newtheorem{example}{Example}
\label{example11}
\example

Solving the Schrodinger equation yields a set of wave functions that satisfy the equation. These wave functions, along with their corresponding eigenvalues, provide information about the quantum properties of the system. The wave functions can be generated using the following formula \cite{Atkins}:

\begin{equation}  \label{ex}
\psi_{n}(x,y) = \frac{1}{\sqrt{n!}} \left(\sqrt{\frac{m\omega}{2\hbar}}\right)^{n} \left(x - \frac{\hbar}{m\omega}\frac{d}{dx}\frac{d}{dy}\right)^{n} \left(\frac{m\omega}{\pi\hbar}\right)^{\frac{1}{4}} e^{\frac{-m\omega x^{2}y^{2}}{2\hbar}}, \quad n=0,1,2,...
\end{equation}

Here, $n$ represents the energy level, and the energy values are given by:

\begin{equation} E_{n} = \left(n + \frac{1}{2}\right)\hbar\omega. \end{equation}

The special case when $n=0$ is called the ground state, with zero energy at the lowest point, and its wave function follows a Gaussian distribution.

A harmonic system, like a particle in a box, demonstrates the discrete energy levels characteristic of the Schrodinger equation. The constants of the equation \ref{ex} in this example are set like this, the particle mass $m=1$, the reduced Planck constant $\hbar=1$, the harmonic oscillator angular frequency $\omega=1$, and the initial potential energy $V_{0}=1$.

During the neural network training, the mean square error is calculated and displayed in a plot, showing the decrease in error over time and iterations. This indicates the neural network's improvement in estimating the wave function. 

We evaluate the performance of the present method and compare it with the Physics-informed neural networks method. The network structure of the proposed method and the $(PINNs)$ method which includes the number of layers, neurons, and activation function are mentioned in the table \ref{arch}


\begin{table}[h]

\begin{centering}
{\scriptsize{}\caption{\label{arch}{ The neural network architecture}
}{\scriptsize\par}}
\par\end{centering}
\begin{centering}
{\scriptsize{}}%
\begin{tabular}{|c|c|c|}
\hline 
 & {\scriptsize{}The present method} & {\scriptsize{}The Physics-informed neural networks method}\tabularnewline
\hline 
\hline 
{\scriptsize{}The number of the inputs} & {\scriptsize{}x and y coordinates} & {\scriptsize{}x and y coordinates}\tabularnewline
\hline 
{\scriptsize{}The number of the hidden layers} & {\scriptsize{}15} & {\scriptsize{}10}\tabularnewline
\hline 
{\scriptsize{}The number of the neurons per of layers} & {\scriptsize{}10} & {\scriptsize{}5}\tabularnewline
\hline 
{\scriptsize{}The number of interactions} & {\scriptsize{}100} & {\scriptsize{}100}\tabularnewline
\hline 
{\scriptsize{}The activation function} & {\scriptsize{}Hermite functions} & {\scriptsize{}Sigmoid function}\tabularnewline
\hline 
\end{tabular}{\scriptsize\par}
\par\end{centering}
\end{table}

\begin{figure}[!h]

\centering 

\includegraphics[width=0.48\textwidth, keepaspectratio]{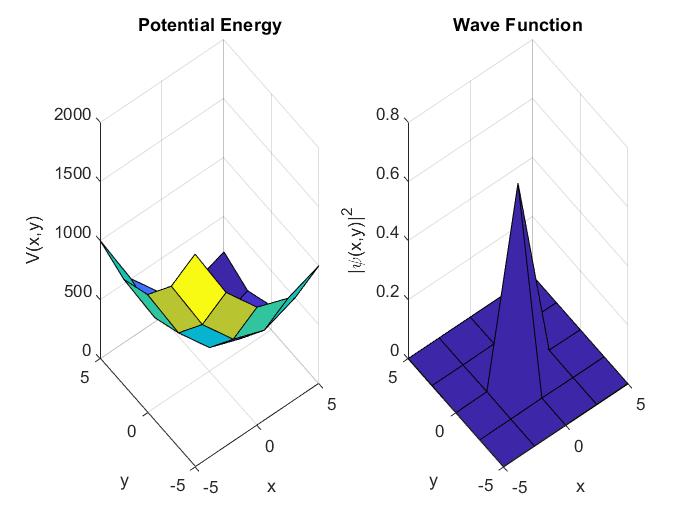}%
\hspace{2mm}%
\includegraphics[width=0.48\textwidth, keepaspectratio]{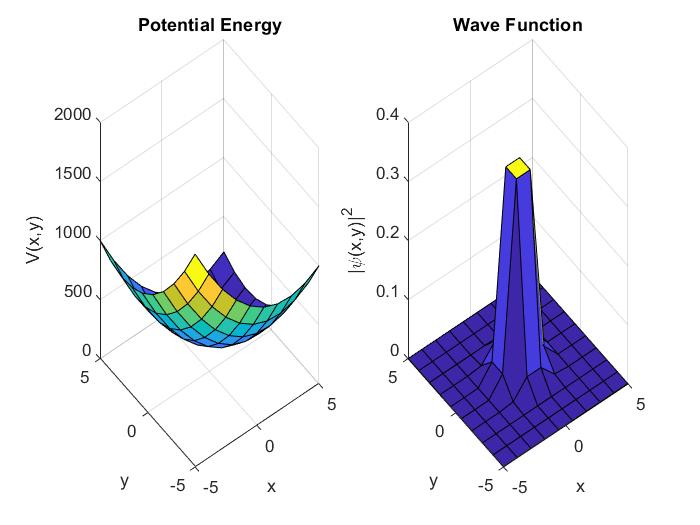}

\vspace{2mm}\vspace{-\lineskip}
\caption{The potential energy plot, and the wave function plot by the Hermite neural network (Left). The potential energy plot and the wave function plot by the PINN (Right)}\label{pewa}
\end{figure}

Another graph compares the actual with the approximated wave function, demonstrating the neural network's estimation proximity to the original wave function.   The plots of the potential energy and wave function obtained from solving the Schrodinger equation using the Hermite neural network and the Physics-informed neural networks are shown in figure \ref{pewa}.

\begin{figure}[!h] 

\centering 

\includegraphics[width=0.48\textwidth, keepaspectratio]{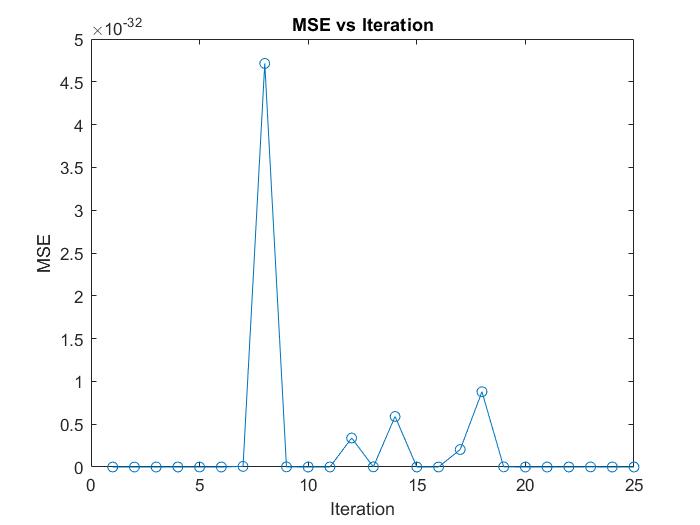}%
\hspace{2mm}%
\includegraphics[width=0.48\textwidth, keepaspectratio]{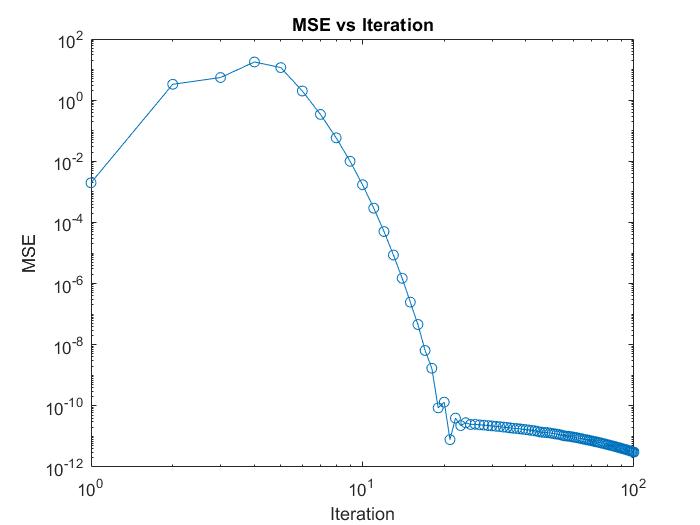}

\vspace{2mm}\vspace{-\lineskip}
\caption{The mean square error by the Hermite neural network (Left). The mean square error by the PINN (Right). }\label{mhp}
\end{figure}

Figure \ref{mhp} illustrates the mean square error from the Hermite neural network and PINN.

\begin{figure}[!h]

\centering 

\includegraphics[width=0.48\textwidth, keepaspectratio]{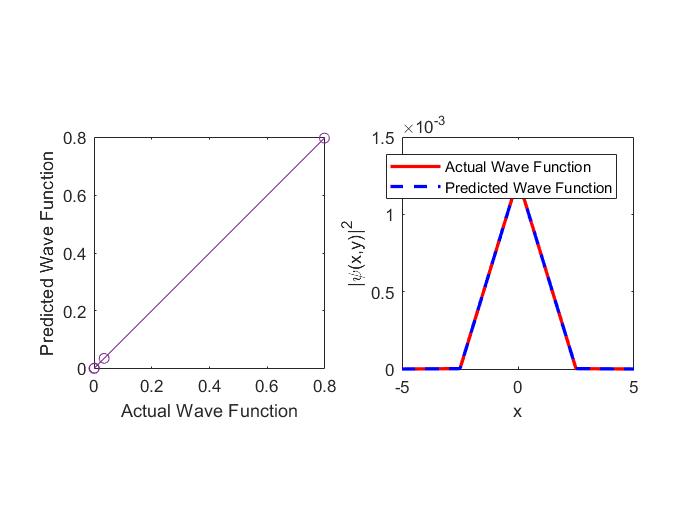}%
\hspace{2mm}%
\includegraphics[width=0.48\textwidth, keepaspectratio]{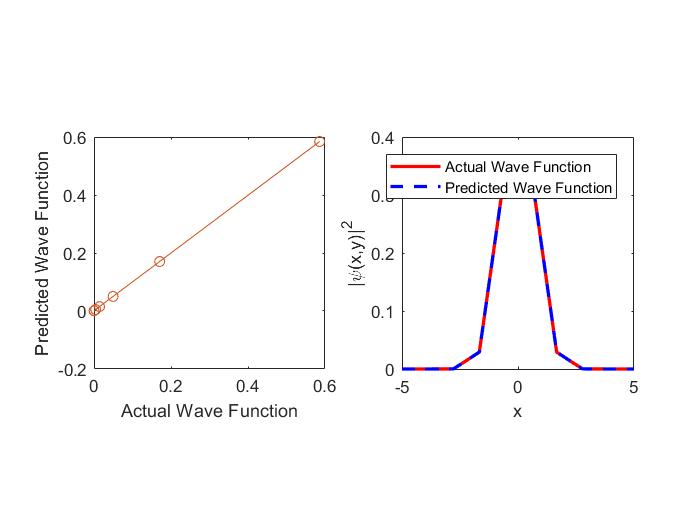}

\vspace{2mm}\vspace{-\lineskip}
\caption{The actual, and the predicted wave function plot by the Hermite neural network (Left). The actual, and the predicted wave function plot by the PINN (Right) }\label{cmh}
\end{figure}

The plots of the approximate and actual wave functions using the Hermite neural network and PINN are in figure  \ref{cmh}.

\begin{figure}[!h]

\begin{centering}
\includegraphics[scale=0.3]{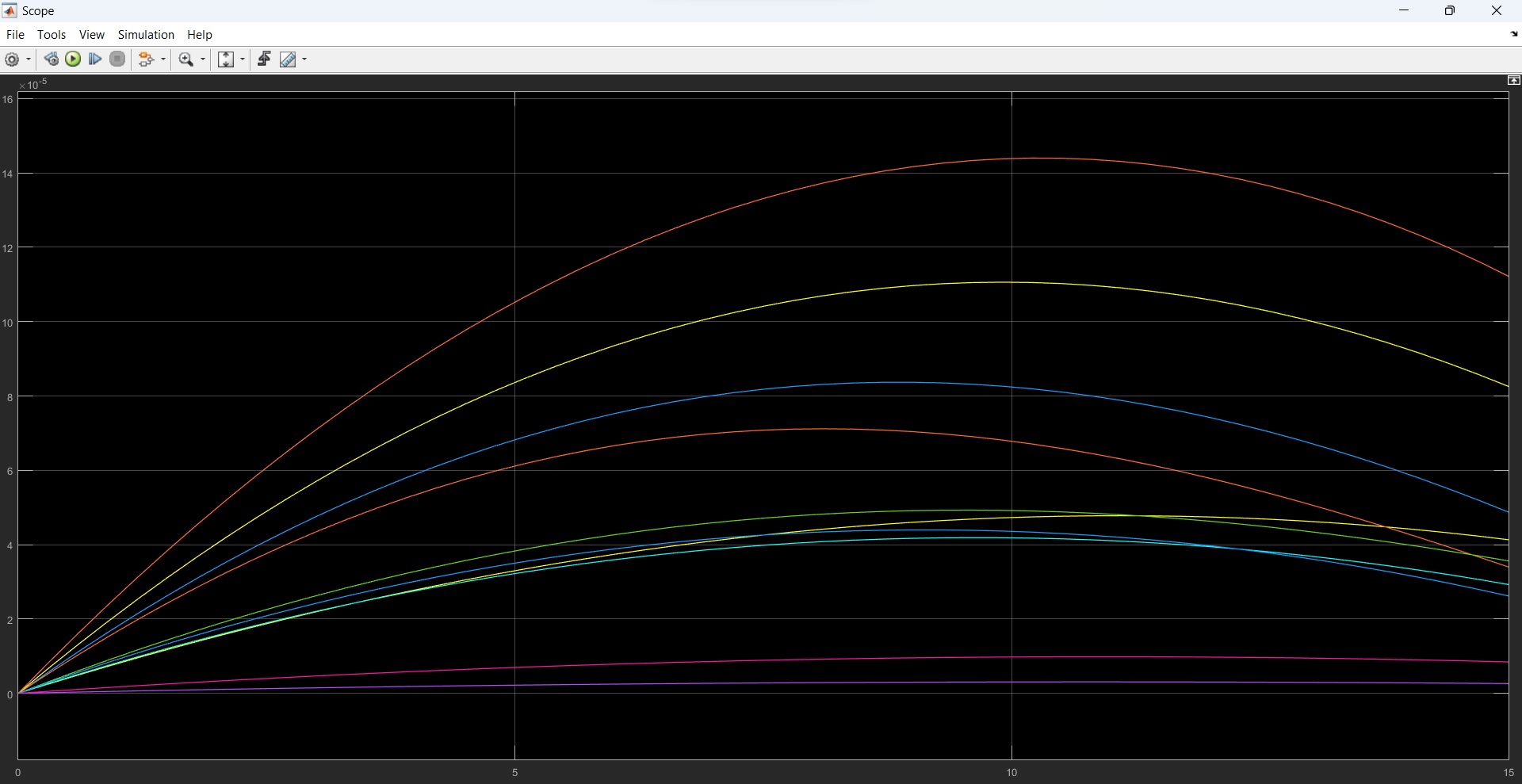}
\par\end{centering}
\begin{centering}
\caption{The Simulink plot by using the present method }\label{simulink}
\par\end{centering}
\end{figure}

In the figure \ref{simulink}, this problem's simulation is shown using the proposed method. 

\begin{figure}[!h]

\begin{centering}
\includegraphics[scale=0.3]{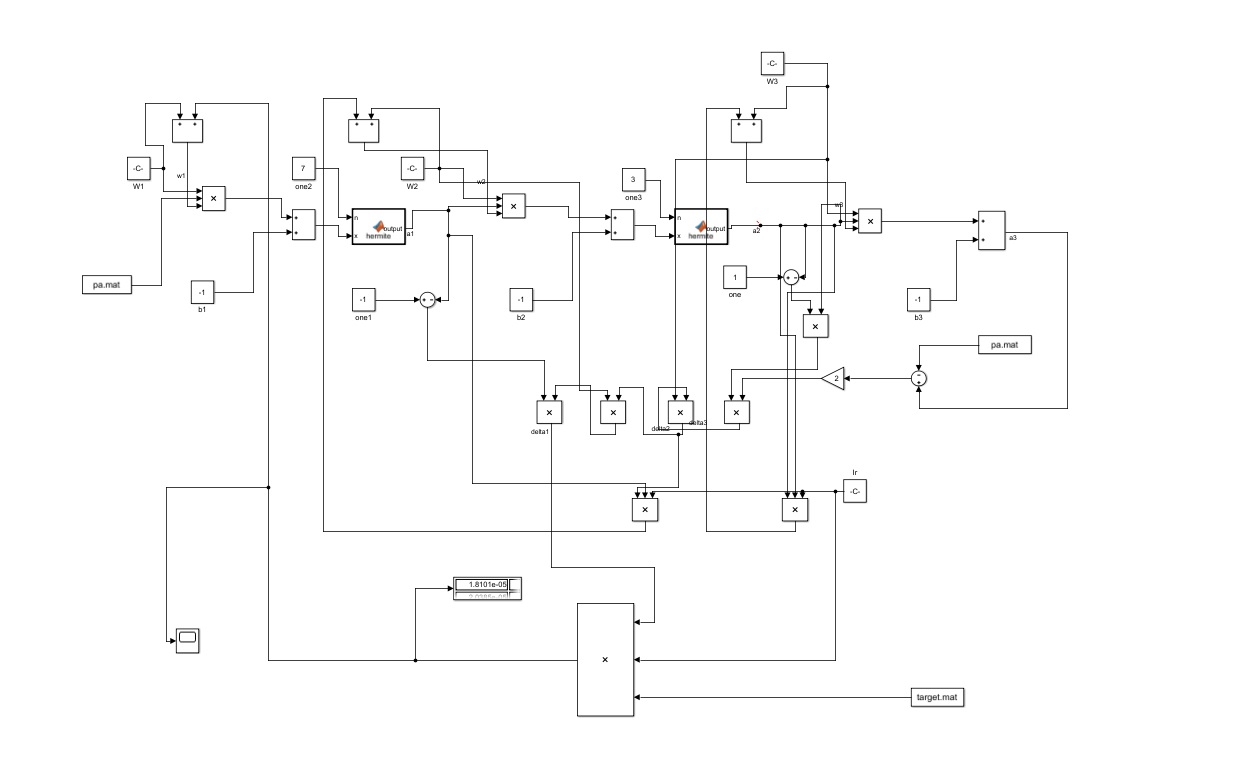}
\par\end{centering}
\begin{centering}
\caption{The Simulink plot by using the present method }\label{herstr}
\par\end{centering}
\end{figure}

This plot \ref{simulink} is for the wave function and shows the approximate location of the electron, also known as atomic orbitals. We showed in figure \ref{herstr} the Hermite neural network structure by the Simulink Matlab.

\vspace{10cm}

\vspace{10cm}

\example

We consider another Schrodinger equation model with the following conditions in this example.  
The complete normalized  wave function is

\begin{equation}
\psi_{n}(x,y)=\frac{2}{L}\sin(\frac{\pi x}{L})\sin(\frac{\pi y}{L}),
\end{equation}

$L$ represents the length of the square coordinates in the Schrodinger equation. $N$ denotes the quantum states corresponding to the variations on each of the coordinates. $(x,y)$ represent the physical coordinates.

This function represents one of the possible states of the wave function in the Schrodinger equation in a 2D box. In this case, you have expressed the general form of the initial wave function $\psi(x,y)$ in this particular state.

The wave function $\psi_{n}(x,y)$ represents the $n-th$ state of a particle in a 2D box. The probability distribution for the occurrence of the particle at points $(x, y)$ is determined using this wave function.

In this wave function, we apply two sine functions at the points$ (x, y)$. The sine function in $x$ is $\sin(\frac{\pi x}{L})$, and the sine function in $y$ is $\sin(\frac{\pi y}{L})$. The coefficient $\frac{2}{L}$ of this function is used for normalization.

This wave function can be used to investigate the probability patterns of the particle in the 2D box. By changing the values of $n$, different patterns of probability distribution in the 2D space can be obtained. Each state represents the points in space where the particle moves with a high probability. Table \ref{arch2} shows the architecture for example 2.


\begin{table}[h]

\begin{centering}
{\scriptsize{}\caption{\label{arch2}{ The neural network architecture}
}{\scriptsize\par}}
\par\end{centering}
\begin{centering}
{\scriptsize{}}%
\begin{tabular}{|c|c|c|}
\hline 
 & {\scriptsize{}The present method} & {\scriptsize{}The Physics-informed neural networks method}\tabularnewline
\hline 
\hline 
{\scriptsize{}The number of the inputs} & {\scriptsize{}x and y coordinates} & {\scriptsize{}x and y coordinates}\tabularnewline
\hline 
{\scriptsize{}The number of the hidden layers} & {\scriptsize{}53} & {\scriptsize{}53}\tabularnewline
\hline 
{\scriptsize{}The number of the neurons per of layers} & {\scriptsize{}15} & {\scriptsize{}18}\tabularnewline
\hline 
{\scriptsize{}The number of interactions} & {\scriptsize{}1000} & {\scriptsize{}1000}\tabularnewline
\hline 
{\scriptsize{}The activation function} & {\scriptsize{}Hermite functions} & {\scriptsize{}Sigmoid function}\tabularnewline
\hline 
\end{tabular}{\scriptsize\par}
\par\end{centering}
\end{table}

\begin{figure}[!h]

\centering 

\includegraphics[width=0.48\textwidth, keepaspectratio]{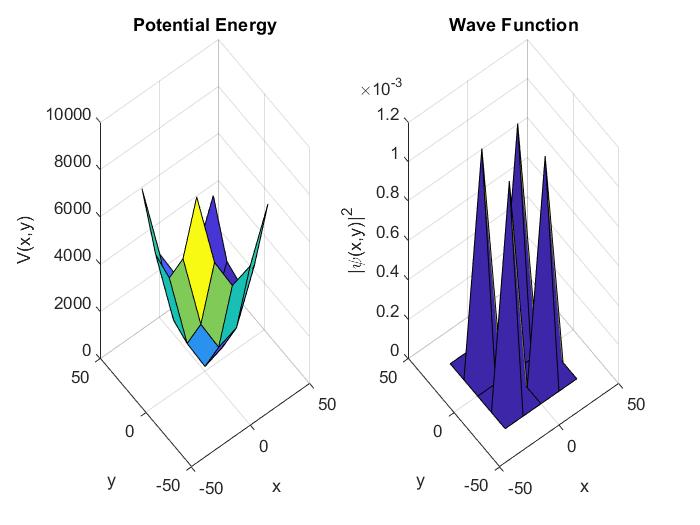}%
\hspace{2mm}%
\includegraphics[width=0.48\textwidth, keepaspectratio]{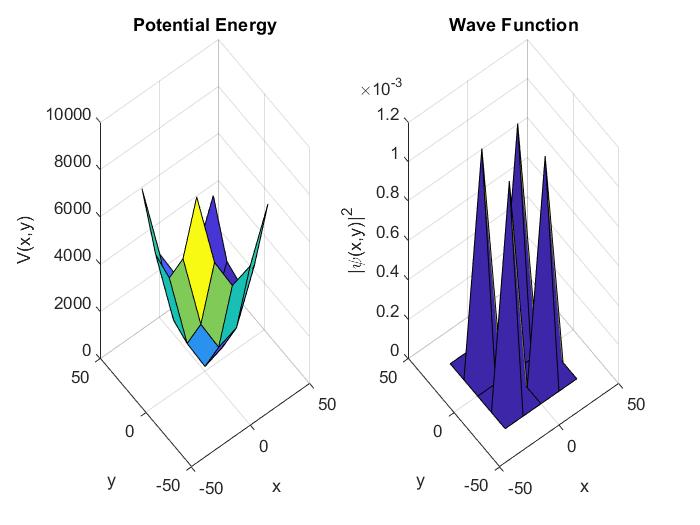}

\vspace{2mm}\vspace{-\lineskip}
\caption{The potential energy plot, and the wave function plot by the Hermite neural network (Left). The potential energy plot and the wave function plot by the PINN (Right)}\label{pess}
\end{figure}

In the potential energy plot \ref{pess}, the horizontal axis represents the quantum numbers indicating different system states, and the vertical axis represents the energy corresponding to each quantum number. The energy-state graph can reveal accepted quantizations, ground energies, and constituent energy components.

\begin{figure}[!h]

\centering 

\includegraphics[width=0.48\textwidth, keepaspectratio]{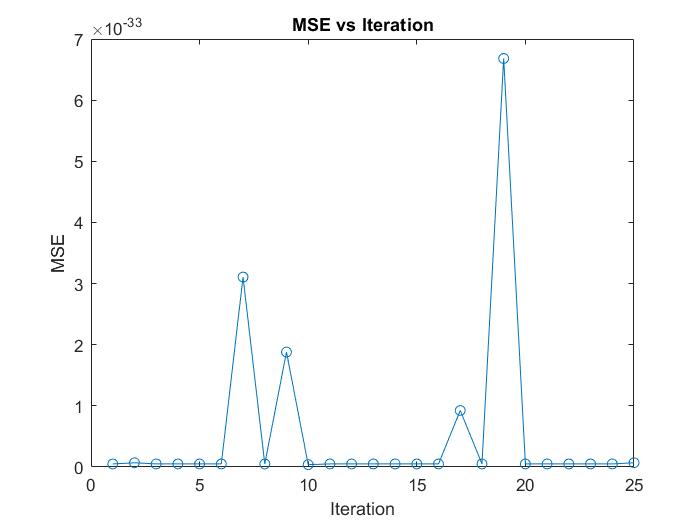}%
\hspace{2mm}%
\includegraphics[width=0.48\textwidth, keepaspectratio]{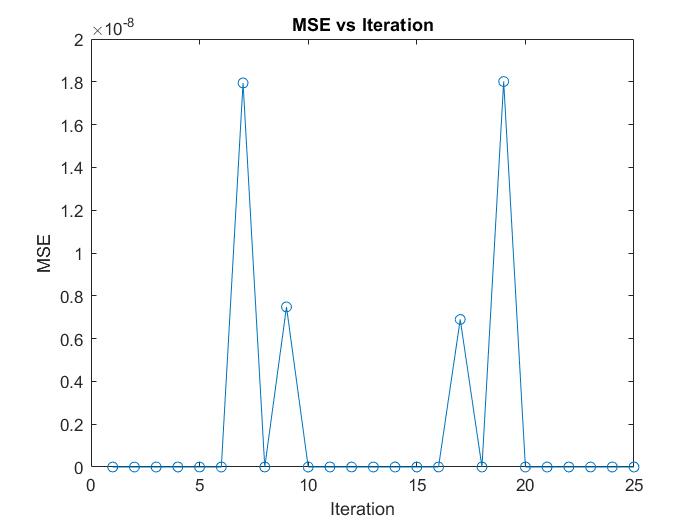}

\vspace{2mm}\vspace{-\lineskip}
\caption{The mean square error by the Hermite neural network (Left). The mean square error by the PINN (Right).}\label{msess}
\end{figure}

The plot \ref{msess} is used to evaluate the performance and training of a neural network. In each iteration of the training algorithm, the corresponding error is computed and recorded in the error history. In this plot, an error matrix is first created with a size equal to the number of training iterations. Then, in each iteration, the corresponding error is computed and recorded in the error matrix. 
The x-axis of this plot represents the number of training iterations, and the y-axis represents the error. Based on the error plot, the neural network is gradually improving in solving the problem. It means that the network is being trained properly.

\begin{figure}[!h]

\centering 

\includegraphics[width=0.48\textwidth, keepaspectratio]{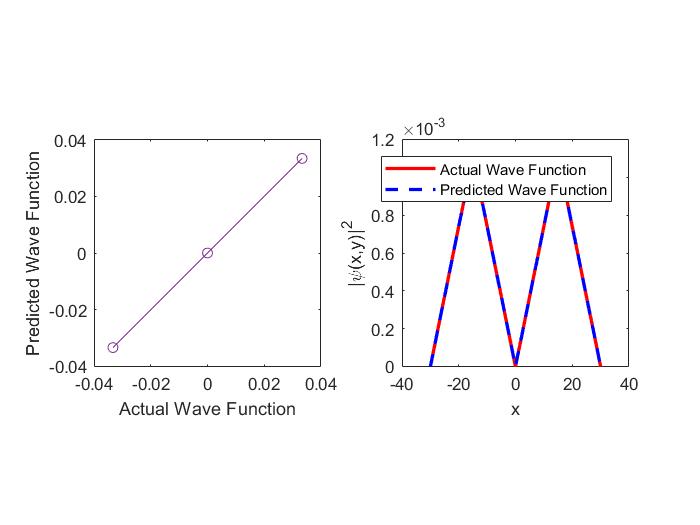}%
\hspace{2mm}%
\includegraphics[width=0.48\textwidth, keepaspectratio]{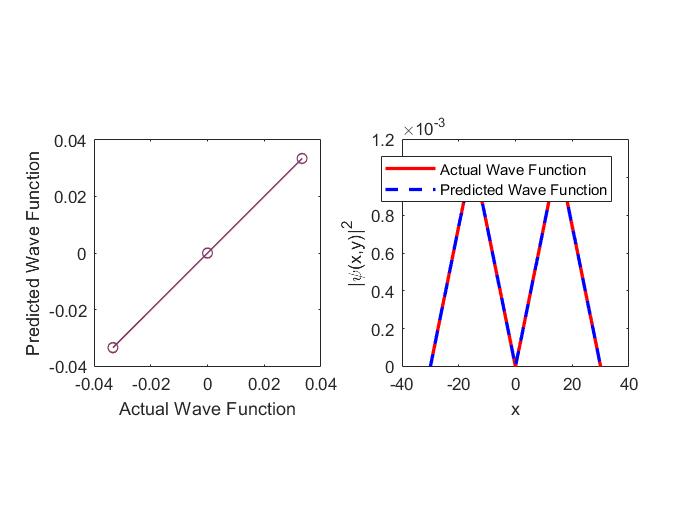}

\vspace{2mm}\vspace{-\lineskip}
\caption{The actual, and the predicted wave function plot by the Hermite neural network (Left). The actual, and the predicted wave function plot by the PINN (Right)}\label{css}
\end{figure}

In Figure \ref{css}, the predicted values of the wave functions by the neural network are compared against the actual values of the wave functions.

\begin{figure}[!h]

\begin{centering}
\includegraphics[scale=0.3]{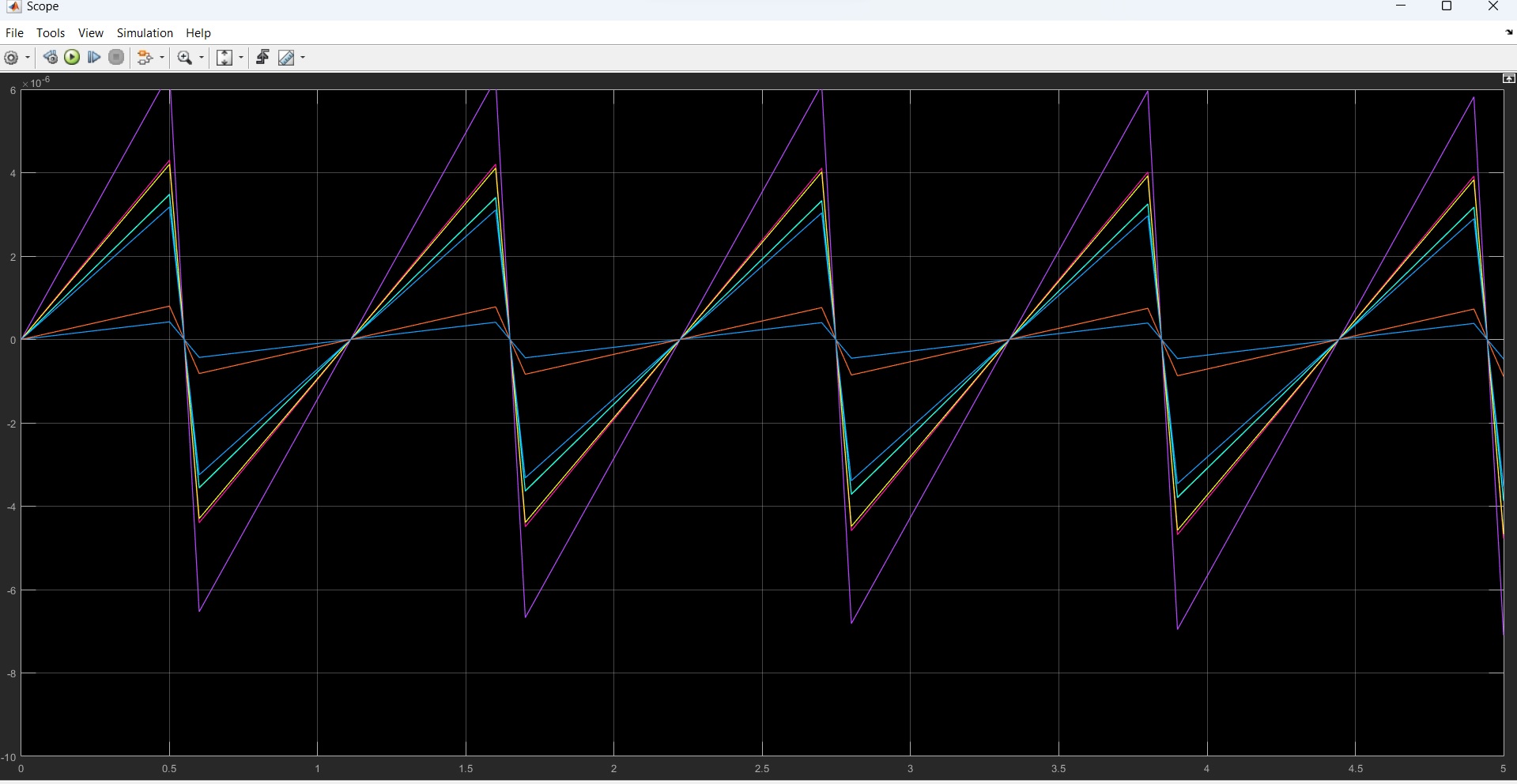}
\par\end{centering}
\begin{centering}
\caption{The Simulink plot by using the present method }\label{shermite}
\par\end{centering}
\end{figure}

Figure \ref{shermite} displays a diagram illustrating the result of a Hermite neural network structure for example 2 in Simulink MATLAB.

\vfill
\vspace{10cm}

\section{Conclusion}

The paper aimed to solve the Schrodinger equation accurately using a mixture of neural networks and the collocation method based on Hermite functions. We have used the roots of Hermite functions as collocation points, which improved the solution's efficiency. Since the Schrodinger equation is defined in an infinite domain, the use of Hermite functions as activation functions led to excellent precision. The proposed method was simulated using MATLAB's Simulink tool. Finally, the results were compared with those from Physics-informed neural networks and the presented method.

\bibliographystyle{unsrt}  
\bibliography{references}

\end{document}